# Probability and Ambiguity

A. Hayrapetyan
(aram.hayrapetyan@gmail.com)

Dec 26, 2024

## Abstract

The title of the article is identical to the title of Chapter 21 in [Gar2001]: because we are going to analyze the probability calculations and the ambiguity of the problem statements. We will analyze 3 out of 4 problems from [Gar2001]: "broken stick", "two boys", "three prisoners"; and the "obtuse random triangle" [Ham1991, p.203-205] problems.

*Keywords*: Broken stick problem, Obtuse triangle problem, Bertrand's problem, Two sons problem, Three prisoners problem, Monty Hall problem; probability function, probability density, randomness, ambiguity

## Introduction

"Charles Sanders Peirce once observed that in no other branch of mathematics is it so easy for experts to blunder as in probability theory" [Gar2001, p.273].

There is an opinion (also expressed in [Gar2001, p.273]) that probability theory fails when the experimental procedure that produces randomness is not precisely defined.

Bertrand's paradox is frequently cited as an example of such ambiguous definitions [Jay1973]. Three typical solutions that produce three different results, are not due to ambiguity of definition or construction method, but rather due to calculating probabilities of different events. Detailed analysis of Bertrand's paradox is available in [Hay2024]: in this article we review some publications and give additional explanations.

The Broken stick problem is also considered as an example of an ambiguous definition: "If a stick is broken at random into three pieces, what is the probability that the pieces can be put together in a triangle? This cannot be answered without additional information about the exact method of breaking to be used" [Gar2001, p.273]. It is the only problem from the 4 reviewed by [Gar2001] that really depends on the breaking method: parallel or sequential. It might be considered a paradox since it is very unintuitive - like the Monty Hall problem. We will analyze the sequential method in detail because, even after correcting [Gar2001, p.278-279] the logic of the original approach, detailed analysis was missing; other authors [Ham1991, p.194-195] , [Ede2015] did not mention sequential method. A simpler model will be used for analysis.

In calculating the probability of random Obtuse triangle creation 2 approaches are analyzed in [Ham1991, p.203-205] and 2 different probabilities obtained. We will review application of the simpler method, used for Broken stick, which yields correct (as in [Por1994] and [Ede2015]) result, and will examine the cause of difference. In addition, we suggest an even simpler method for calculating the Obtuse triangle probability.

## Randomness

To eliminate ambiguity in the definition of randomness and prevent philosophical discussions on its undefinability [Cha2006, pp. 123-125, Ham1991, pp.20-22], we will use Laplace's definition of probability of a random event: "the probability of an event is the ratio of the number of favourable cases to the number of all possible cases, when nothing engenders a belief that any one of these cases should occur rather than any other, which renders them, for us, equally possible" [Mar1994] (see also the *Symmetry as the Measure of Probability* section in [Ham1991, pp.9-10]).

We adhere to this definition of randomness for practical reasons. By analysing ambiguities we want to understand the challenges of probabilistic approach to experimental testing of theoretical ideas and identifying possible pitfalls. "Randomness can never be proved, only the lack of it can be shown" [Ham1991, p.20].

# Bertrand's problem

## Literature Review

This review is not exhaustive or thorough: apologies to uncited authors, because below we make tacit assumption that their analysis is similar to the cited, and to the reviewed authors, because we might have misread or misinterpreted their ideas (there are too many long expressions containing integrals in the papers).

The majority of authors explain paradox by ambiguity of the problem definition [e.g. Gar2001, Jay1973, Mar1994] or, what seems equivalent, by not specifying the method of generating random events. "Each of these procedures is a legitimate method of obtaining a "random chord." The problem as originally stated, therefore, is ambiguous. It has no answer until the meaning of "draw a chord at random" is made precise by a description of the procedure to be followed" [Gar2001, p.277].

Others point to semantic issues [Mar1994, DiP2010, Pet2019, Che2023]: what is random? Some of these [Mar1994, DiP2010, Pet2019] see the problem in "(geometrical) probabilities defined on uncountably infinite sets" [Pet2019].

To explain or resolve the "paradox" concepts like principle of indifference and its extension, the principle of maximum entropy [Wan2011], principle of "maximum lack of knowledge"and Laplace's principle of insufficient reason (or of indifference) [Aer2014] are discussed.

In [Ham1991, p.201] an interesting Bertrand's reversed problem is analysed. However, it is a reverse of just one method: instead of rolling a broom handle towards a circle [Gar2001, p.276], a coin is tossed on a line.

In Bertrand's case the favorable events are chords longer than the side of the inscribed equilateral triangle and all possible cases - sample space - are all possible chords (and not just perpendicular to the radius - see [Hay2024]).

There should be no concerns about the inapplicability of Laplacian probability definition to continuous cases because we know from experience (or measure theory) that the probability of a randomly tossed point at a line of length L to land on a section of length l of the line is l/L. The key is to calculate ratios for the same type objects, of same type measures: length vs. length, surface vs. surface, etc. "In the infinitesimal notation of scientists and engineers we can talk about the probability of being in an interval dx as being p(x) dx, but not about the probability at a point. Talking about the probability density p(x) is very convenient, but careful thinking requires us to go to the integral form whenever there is any doubt about what is going on" [Ham1991, p.202].

No author explains why the problem is ambiguous or not well-posed (with some exception, e.g. [Gar2001, p.273]). Blaming undefined processes is not justified, because the 3 processes specified by Bertrand produce different events, different types of chords, or different sample spaces. The chord, equilateral triangle, circle, length more than another length, and random (rather non-random) are perfectly defined mathematical notions. Instead of specifying what is exactly ambiguous in the statement "We randomly draw a chord in a circle" [Ber1889, p. 4], the researches started throwing broomsticks, straws, pebbles towards a circle or coins at lines (reverse problem) and calculating probabilities of events that have nothing (or not much) to do with what Bertrand was asking - draw random chords. Unfortunately, Bertrand himself is the source of confusion.

One can consider a chord perpendicular to the radius at random point on a circle (or radius) special, non-random (which makes sense: it emphasizes lack of randomness [Ham1991, p.20]) or consider all three special types generated randomly. These assumptions can be classified as semantic issues, as matters of definition. That is why they are not important. What is important is that Betrand specified three different sample spaces.

Bertrand's problem has been used in a study on human behavior (?).

"Apparently nothing resembling any of the three procedures is actually adopted by most people when they are asked to draw a random chord. In an interesting unpublished paper entitled "The Human Organism as a Random Mechanism" Oliver L. Lacey, professor of psychology at the University of Alabama, reports on a test which showed the probability to be much better than 1/2 that a subject would draw a chord longer than the side of the inscribed triangle" [Gar2001, p.277].

This is the Abstract of the mentioned above article: "98 Ss were instructed to draw a chord (or set of chords) at random across a circle with an inscribed equilateral triangle. The chords were measured, and the principal result was the finding that the proportion of chords found to be greater than the length of the side of the inscribed equilateral triangle was .628 — significantly greater than the upper predicated mathematical limit of .5. Perceptual effects were suggested to account for the finding" [Lac1962].

Why was the circle with the inscribed triangle chosen? A much more interesting experiment is asking the subjects to draw a chord in the empty circle (also a "chord" in other shapes: square, triangle) without giving any extra information. Such tests will not contribute to understanding Bertrand's problem, but they might (will) contribute to understanding human behavior.

### *Explanation*

The chord constructed using second and third methods is not random, but rather special (which is the opposite of being random). It is built at a randomly chosen point, but its direction is uniquely defined by that point. In other words, only one of the two degrees of freedom is random. Thus defined - perpendicular to the radius - chord is special: it is always larger than the side of the inscribed equilateral triangle, if the point is inside the circle of radius r/2, and always smaller - when the point is outside the radius r/2. That is why the probability for the second method is equal to the probability of randomly selecting a point inside an r/2 circle and using the third method - to the probability of randomly selecting a point on the half of a straight line segment.

There is no need for additional information. It is a perfectly formulated problem with a uniquely identifiable probability of 1/3. The other 2 probabilities can be interpreted as probabilities of different objects or as a probability of different sample spaces. The 3 methods calculate probabilities of these 3 objects: a random chord (Method 1), a chord perpendicular to the radius in a randomly selected point on the radius(Method 2), a chord perpendicular to the radius in a randomly selected point in the circle (Method 3).

## *Broken stick problem*

### *Selecting 2 points at once*

The Broken stick problem is considered a classic example of an ambiguous definition [Gar2001, p.273]. If you select 2 random points at once, then the probability is 1/4; if sequentially, it is 0.193.

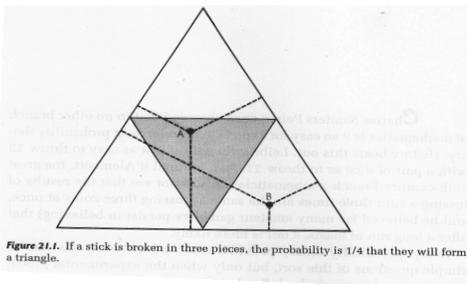

Figure 21.1. If a stick is broken in three pieces, the probability is 1/4 that they will form a triangle.

"One method is to select, independently and at random, two points from the points that range uniformly along the stick, then break the stick at these two points. If this is the procedure to be followed, the answer is 1/4, and there is a neat way of demonstrating it with a geometrical diagram. We draw an equilateral triangle, then connect the midpoints of the sides to form a smaller shaded equilateral triangle in the center (see Figure 21.1). If we take any point in the large triangle and draw perpendiculars to the three sides, the sum of these three lines will be constant and equal to the altitude of the large triangle. When this point, like point A, is inside the shaded triangle, no one of the three perpendiculars will be longer than the sum of the other two" [Gar2001, p.273].

A simpler model (diagram) for this approach is a line **o**=1, which represents the perimeter of the triangle. To build a triangle from the pieces (sections) of the line **o** we need to ensure that each side is less than 1/2. There are three configurations, when 2 cuts make no triangle: 1) both cuts are in the first half of the perimeter length stick, 2) both cuts are in the second half of the stick, and 3) the distance between cuts is longer than 1/2. Thus the probability for *not* getting a triangle is 3/4, which makes the requested probability 1/4.

### *Selecting 2 points sequentially*

If we break the stick once and pick the smaller piece, we cannot build a triangle because the sum of two sides will be smaller than the third side. The probability of selecting the bigger part for the second break is 1/2.

After selecting the larger piece of the stick, we need to calculate the conditional probability for making the requested cut $c_2$, considering the outcome of the first cut $c_1$:

1. The range for x is [0-:-1/2] (for calculations simplicity we assume that the length of original stick is 1)
2. The range for making the $c_2$ cut is [0-:-(1-**x**)]

3. These are the lengths of 3 edges of triangle after cut $c_2$: a) **x**, b) $c_2$, and c) 1-**x**-$c_2$

For a continuous variable, we need to identify the probability function to calculate probability density. The probability function is the curve that represents probability values of the variables, while the probability density is the integral over the interval of the variable values, analogous to the probability distribution for discrete variables.

We will use Diagram 1 - which is a redraw of Figure 21.1 above - to calculate the probability function for the second cut in sequential case.

First couple of geometrical details about the equilateral triangle and the relevance to the problem. The stick length equals the height IC of the big triangle on Diagram 1. We assume it is equal to 1. The sum of lengths of perpendiculars to the sides from any point inside the big triangle - from points like A and B in Figure 21.1 or F in Diagram 1 - equals the height of the triangle.

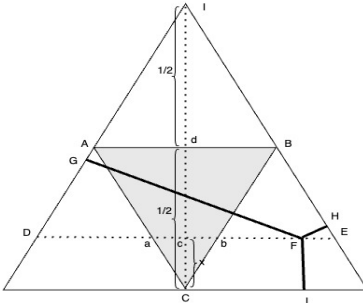

Diagram 1

Any side *s* of equilateral triangle equals **h** · 2/√3.

Let's assume that the length of the first cut is **x** and it is the smaller part of the two: less than half of the stick as it is shown on Diagram 1. After the first cut the choice of point for the second cut is limited to the points on the DE. In the original approach [Gar2001, p.273-p.275] it was assumed that the point can be selected in any location on the surface below the AB line. Hence, the [incorrect] probability of 1/2 · 1/3 = 1/6 is calculated.

To satisfy the requirement of the problem the second point - as it was in the case of simultaneous cuts before - should be inside the shaded triangle. However, now we have an additional condition to meet - the point should be on the DE line.

The probability for that is the ratio of ab/DE. These lengths expressed in term of **x** are:

$$s_{ab} = 2x/\sqrt{3} \quad [1]$$

$$s_{DE} = 2(1-x)/\sqrt{3} \quad [2]$$

The probability function for getting requested second cut $f_p$ of the length **x** is:

$$f_p = x/(1-x) \quad [3]$$

The probability density for the [3] is

$$\int_0^{1/2} (x/(x-1))\, dx = (-x - \ln|1-x|) = -1/2 - \ln 1/2 = -0.5 + 0.693147 = 0.193147 \quad [4]$$

In the corrected analysis of the problem [Gar2001, p.279] the probability function is defined using *log*, which is confusing because from the context (that does not have the explanations above) it is not obvious that $log_e = ln$ is meant.

We can use a stick **o** of length 1 (see in previous section) to derive the expression [3]. After the first cut we have 2 sticks: short stick of length **x** and long - 1-**x**. The condition of a successful second cut is to get it in the range [(1/2-x), (1-x)-(1/2-x)] on the 1-**x** long line section. This section is (1-x)-(1/2-x) - (1/2-x) = **x** long.

*Explanation*

The fact that we get different probabilities in case of sequential and simultaneous cuts feels paradoxal. However, the opposite might seem paradoxical to other people. For example, "Jean le Rond d'Alembert, the great 18th-century French mathematician, could not see that the results of tossing a coin three times are the same as tossing three coins at once" [Gar2001, p.273].

Then why for coins the probability for sequential and simultaneous tosses are not different, while for breaking a stick they are? In coin tosses all three events are independent. In the case of breaking a stick the simultaneous cuts are independent. The probabilities of both favorable cuts are the same; in sequential case the probability of a favorable second cut depends on the length **x** of the first cut. It is dependent, conditional probability.

## Obtuse random triangle problem

### Overview

"Three points are taken at random on an infinite Plane. Find the chance of their being the vertices of an obtuse-angled Triangle" [Dod1893, p.14].

In [Ham1991] 2 methods of calculating the probability of drawing an obtuse triangle on a plane are analyzed: 1) by placing "the x-axis of a cartesian coordinate system along the longest side of the triangle, and the origin at the vertex joining the second longest side. We will scale the longest side to be of unit length since it is only the ratio of areas that will matter" [Ham1991, p.203], 2) by assuming "the x-axis along the second longest side with unit length, and the vertex at the junction of that side with the shortest side. [Ham1991, p.204]. Fixing the length of the small side makes probability calculations impossible because the lengths of long and medium sides change in the [1, ∞] range.

We will use elementary geometry instead of calculus (integrals) and will not align any side with x-axis since this does not seem important and it does not simplify calculations. The coordinate axis and the triangle exemplar are presented on diagrams for convenience. On the Diagrams 2 and 3 below we will denote short, medium (second long), and long sides **s**, **m**, and **l** and call calculation methods M and L respectively.

### L method

On the Diagram 2 we have an example of a random ABC triangle in the Cartesian coordinate plane. We assume that the length of the long **l** (AC) side is equal to 1. If the vertex B of the triangle is located inside the half circle AC of radius **r**=AE=EC=1/2, then the triangle is obtuse. Per assignment of sides of the Diagram 2 the point B cannot go beyond the line DE: the inequality **s** <= **m** will be violated. Thus, the area in question is a half of a half circle AC or quarter of the full circle. For the same reason B can go above the radius **r** =**l**/2 circle, but not to the right of DE, even though the figure ACD, limited by the arcs AD and CD of radius **R**=AC=CD=1, is the area of all possible triangles. Due to symmetry the peculiarity of side assignments on the Diagram does not affect the end result: the probability of getting an obtuse triangle is equal to the ratio of the area of the half circle of radius **r** to the area of figure ADC or area $A_o$ of quarter circle of radius **r** to the area $A_A$ of figure ADE - area of all possible triangles) .

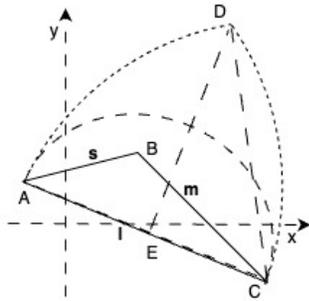

The area $A_o$ of obtuse triangles is π/16, while the area $A_A$ is the area of figure ADE, which is 1/2(π/3 - √3/4): the area of the section ACD of the R=1 radius circle minus the area of triangle CDE

Diagram 2       $P_o = A_o / A_A = π/16 / (1/2 \cdot (π/3 - √3/4)) ≈ 0.639$

The value $P_o = 3 / (8 - 6/π √3) ≈ 0.64$ is obtained originally in [Dod1893, p.83-84] by using the Diagram 2.

### M methode

Frank Wattenberg suggested this solution in 1973 [Por1994].

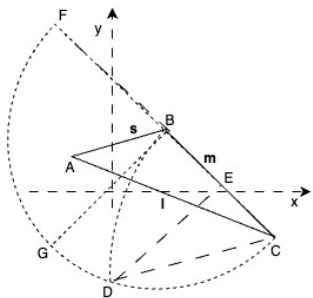

On the Diagram 3 we consider the side **m** is equal to 1. In this case the triangle is obtuse if the vertex A is in section BFG, which is a quarter of the circle of radius r=BC=BF=BG=1. It means that $A_o$ is π/4. The long side **l**=AC cannot go to the right of the arc BD or outside the half circle CDGF, because it cannot be shorter than **m**, while **s** cannot be longer than **m**. It means that the area $A_A$ is equal to the area of figure BDGF.

It is not hard to notice that the figure BCD of Diagram 3 is identical to the figure ACD of Diagram 2, which means that the area of the figure BDG is equal to π/4 - (π/3 - √3/4). Thus $A_A$ = π/4 - (π/3 - √3/4) + π/4 = π/6 + √3/4.

$$P_o = A_o / A_A = π/4 / (π/6 + √3/4) ≈ 0.821…$$

Diagram 3       Note. The value π/6 - √3/2 obtained in [Ham1991, p.204] for $A_A$ is incorrect (due to error in integral construction earlier). However, the resulting value of $P_o$ is [magically] evaluated correctly.

*Simple model*

It is evident and easy to prove that scale does not affect the results [Ham1991, p.204] - choice of unit plays no role in making a triangle by cutting a stick or in making a triangle obtuse. We care only about the ratios of sides (angles).

We consider a straight line segment named π representing 180° and randomly place two points on it. What is the probability that one of the 3 resulting sections is bigger than half of the π? For 2 points to land on the first half of the π is 1/4. Same is true for the second half. The section between the 2 tossed points can also be longer than π/2. Hence, the requested probability is 3/4.

This model is identical to the Broken stick problem, where a unit line represents the sum of side lengths (perimeter), while in this case, it represents the sum of angles

The 3/4 result is obtained and proved in [Por1994] analytically by considering selection of 6 points in $R^6$ space and using "four well-known facts about normal sampling": the probability of an angle in a triangle to be ≥ 90° is 1/4, hence the requested probability is 3/4. Using Figure 21.1 (or Diagram 1) - see in the Broken stick problem section - the 3/4 is calculated [Ede2015]. If we consider the big triangle of Figure 21.1 representing the area of all triangles, then a small gray triangle corresponds to acute triangles, while the other 3 small white triangles - to obtuse.

*Simpler model*

A triangle can have only one big angle and only this angle can be obtuse. The range of uniform distribution of the big angle value is [π/3, π]. This is the sample space. Favorable outcomes are in the [π/2, π] range. Hence, the probability of randomly getting an obtuse triangle is 1/2 ⋅ π / 2/3 ⋅ π.

*Explanation*

Probabilistically, the Broken stick and Obtuse triangle problems are equivalent. The sequential and simultaneous methods are applicable to either of them. We get 3/4 per simple method, which is different from L and M methods.

Why are we getting different results? Which is the correct one? "What is wrong? It is not the choice of coordinate systems, it is the idea that you can choose a random point in the plane. Where is such a point? It is farther out than any number you can mention! Hence the distance between the points is arbitrarily large -- the triangle is practically infinite and the assumptions of the various steps, including supposing that the side chosen is of unit length, are dubious to say the least" [Ham1991, p.205].

Assigning unit length to one side defines different ranges for the lengths of the other two sides:

1. L method: l = 1, then s = [0 : 1], m = [0 : 1],
2. M method: m = 1, then s = [0 : 1], l = [1 : 2],

These methods define different [shapes of] sample spaces (see Diagrams 2 and 3), which leads to different probabilities for the same problem.

It is significant that the L method "was accepted for so long, and it indicates the difficulty of using probability in practice, especially continuous probability. You need to be very careful whenever you come in contact with the infinite." [Ham1991, p.204].

How can one model the probability of selecting a random value on an infinite line? "Clearly, the problem is not well-posed since the notion of "at random on an infinite plane" is not defined precisely. … it seems reasonable that one consequence of taking points "at random in the plane" is that the induced distribution in $R^6$ be spherically symmetric" [Por1994]. Using transformation groups rather than spherical symmetry in $R^6$ yet another result is obtained - 0.8426, but building a triangle by randomly tossing 3 points on the circle still gives 3/4 [Por1994]. "None of these alternative approaches uses a measure in $R^6$. Thus, the connection with spherical symmetry in $R^6$ is unclear, and the fact that the answer 3/4 arises so often seems to be fortuitous. It would be extremely interesting to find a general principle underlying all the P = 3/4 answers. Note also that none of these approaches deals with choosing a random element from the set of triangles. The sample space is always strictly larger than the set of triangles. In fact, it is rather difficult (or perhaps impossible) to find a natural group operating transitively on the set of triangles (however this set is defined)" [Por1994].

We cannot model uniform distribution on an infinite plane neither physically, nor computationally because "physical reality" and computers set boundaries to sample space. For example, in a computer model polar coordinate radius or cartesian axes are limited by some max value specific to the operating system. On a computer, you end up modeling uniform distribution on a plane represented either by circle or square (or other figure).

Simulation on a circle or a square yields different results (see Appendix). "If the difference between the answers for the square and circle is not obvious, then consider a long, flat rectangle and let the size approach infinity while keeping the "shape" of the rectangle constant. The probability of an obtuse triangle in a long, low rectangle is very high and depends on the shape of the rectangle. Hence, just which shape should we choose to let approach infinity? The answer clearly depends on this choice of shape and not on the size of the figure" [Ham1991, p.205].

*Experimentation*

To investigate the Obtuse triangle problem, we implemented multiple computer routines to model the sample space - see in Appendix.

We included the statistics for the *Simpler model* (labeled as *Generated* in the Appendix Table) to compare the distributions of **s/m** and **m/l** side length ratios (SdM and MdL respectively in the Table) with other triangle creation models. The results in the Table match respective values of several calculations made before - see, for example, [Ban2009] (except for normal distribution).

## Two boys problem

"Readers were told that Mr. Smith had two children, at least one of whom was a boy, and were asked to calculate the probability that both were boys. Many readers correctly pointed out that the answer depends on the procedure by which the information "at least one is a boy" is obtained. If from all families with two children, at least one of whom is a boy, a family is chosen at random, then the answer is 1/3. But there is another procedure that leads to exactly the same statement of the problem. From families with two children, one family is selected at random. If both children are boys, the informant says "at least one is a boy. If both are girls, he says "at least one is a girl." And if both sexes are represented, he picks a child at random and says "at least one is a ...," naming the child picked. When this procedure is followed, the probability that both children are of the same sex is clearly 1/2. (This is easy to see because the informant makes a statement in each of the four cases - BB, BG, GB, GG - and in half of these cases both children are of the same sex.)" [Gar2001, p.277].

For families with two children, including at least one boy (like the Smiths), the probability of both children being boys is 1/3, if the requested event is the *two boys from all families with at least one boy*. In this case the GG is removed from the sample space [Ham1991, p.25]. If the request is: What is the probability in a *family with 2 children of having 2 boys when one of them is a boy* ? Then the probability is 1/2, because the probability that a family with 2 children has 2 boys is 1/4.

The ambiguity in this problem arises from the definition of the sample space. If the task is to calculate the *unconditional* probability of having 2 boys among the families with at least one boy, then the answer is 1/3. But if the question is: "What is the *conditional* probability of having 2 boys when one of them is a boy among all families with 2 children?", then the answer is 1/2.

## Three prisoners problem

"Three men - A, B, and C - were in separate cells under sentence of death when the governor decided to pardon one of them. He wrote their names on three slips of paper, shook the slips in a hat, drew out one of them, and telephoned the warden, requesting that the name of the lucky man be kept secret for several days. Rumor of this reached prisoner A. When the warden made his morning rounds, A tried to persuade the warden to tell him who had been pardoned. The warden refused.

"Then tell me," said A, "the name of one of the others who will be executed. If B is to be pardoned, give me C's name. If C is to be pardoned, give me B's name. And if I'm to be pardoned, flip a coin to decide whether to name B or C."

"But if you see me flip the coin," replied the wary warden, "you'll know that you're the one pardoned. And if you see that I don't flip a coin, you'll know it's either you or the person I don't name."

"Then don't tell me now," said A. "Tell me tomorrow morning."

The warden, who knew nothing about probability theory, thought it over that night and decided that if he followed the procedure suggested by A, it would give A no help whatever in estimating his survival chances. So next morning he told A that B was going to be executed.

After the warden left, A smiled to himself at the warden's stupidity. There were now only two equally probable elements in what mathematicians like to call the "sample space" of the problem. Either C would be pardoned or himself, so by all the laws of conditional probability, his chances of survival had gone up from 1/3 to 1/2" [Gar2001, p.277-278].

This is similar to a well known and well analyzed Monty Hall problem, but the description of the problem is more complicated. The probability of a pardon for prisoner A remains 1/3, however, the probability for prisoner C it is 2/3, because before gaining the information about the sad fate of prisoner B the probability of B or C pardoned is 2/3, and it remains the same after gaining information. But since we know that the probability of B being pardoned is 0 the "whole" 2/3 go to prisoner C. This situation would have been absolutely identical to the Monty Hall problem, if the governor's order would have been pardoning the person occupying cells A, B, or C - writing down cell numbers, rather than prisoners names. In this case if prisoner A could have managed to convince the warden and prisoner C to switch cells, then the probability for her to be pardoned would have risen twofold - 2/3.

In both cases the warden and the show host knew for sure the fate of the prisoner or contestant. In the [Ham1991, p.31-32] the game show problem is analyzed also for the case when the host acts probabilistically. The game show example is relatively simple for analysis. However, for more difficult problems the "patient reasoning we are using will not work very well" [Ham1991, p.32].

## *Conclusions*

It is evident that probability problems often have counterintuitive, difficult-to-grasp solutions:

"When I introduced the three prisoners paradox in my October 1959 column, I received a raft of letters from mathematicians who believed my solution was invalid. The number of such letters, however, was small compared to the thousands of letters Marilyn vos Savant received when she gave a version of the problem in her popular Parade column for September 9, 1990" [Gar2001, p.283]. "The red-faced mathematicians, who were later forced to confess they were wrong, were in good company. Paul Erdös, one of the world's greatest mathematicians, was among those unable to believe that switching doors doubled the probability of success. Two recent biographies of the late Erdös reveal that he could not accept Marilyn's analysis until his friend Ron Graham, of Bell Labs, patiently explained it to him" [Gar2001, p.284].

The 3 different outcomes in *Bertrand's problem* correspond to 3 different sample spaces. We can only guess if by stating that "None of the three is false, none is correct, either can be questioned" [Ber1889, p. 5] Bertrand meant to show that probability depends on the procedure (like in the case of *Broken stick*) or that geometric probabilities are problematic. If he did, then - it seems - he failed.

The *Three prisoners* (or *Monty Hall*) problem probably has the most unintuitive solution or explanation.

The definition of sample space in the above formulation of the *Two sons* problem is ambiguous.

The *Broken stick* problem really depends on the breaking method. If there is any physical "random" process that "aims" to cut a stick in three pieces and build a triangle then the chances are 0.193… . However, this is only true if the process operates on only one piece at a time: like an insect cutting a leaf stem. The odds for a meteorite shower seem much better. The issue is not in simultaneous vs. sequential, but rather in sequential (dependent) vs parallel (independent) processing. This explanation - though satisfactory - still remains very unintuitive.

There is another issue concerning the sequential method. How realistic is it? Is the sequential method relevant to any natural, real process? [This might be the reason for ignoring analysis by some authors.]

The above 4 problems can be considered as resolved and explained. They can be tested (simulated) physically and computationally. The *Obtuse random triangle* problem in its original formulation can be not.

The *Broken stick* and the *Obtuse random triangle* problems are equivalent when looked at as a triplet of real variables of constant sum. This method is sensitive to the processing approach: parallel or sequential. The results for either can be tested and computer modeling confirms theoretical calculations.

Is it possible to simulate uniform distribution on an infinite line, plane, etc.? To avoid singularities we can identify variables for which relations, for example ratios, rather than absolute values are important, and they do not impose "boundaries" on an "infinite plane". Since different approaches using unit length stick or sum of angles to limit variable ranges give 3/4 then it is the right probability for the *Obtuse triangle* (and 1/4 for the *Broken stick*). This solution does not seem ambiguous, unintuitive, paradoxical, or "equally good (bad)" [Por1994]. It is just impossible to test, confirm experimentally.

## *Appendix - Random triangle creation*

The results of modeling different random distributions on a plane (polar coordinates) by Java *Random* class:

| Method \ Metrics | $\mu$ | Md | Min | Max | $\sigma^2$ | Skew |
|---|---|---|---|---|---|---|
| **Generated** | P = 0.75 | | | | | |
| SdM | 0.49 | 0.47 | 0 | 1 | 0.09 | 5.17 |
| MdL | 0.83 | 0.86 | 0.5 | 1 | 0.02 | -32.02 |
| ρ | 0.67 | 0.71 | 0 | 1 | 0.06 | -31 |
| Θ | 0 | 0.01 | -3.14 | 3.14 | 3.29 | -0.21 |
| **Bertrand #2** | P = 0.7573 | | | | | |
| SdM | 0.58 | 0.59 | 0 | 1 | 0.07 | -10.5 |
| MdL | 0.78 | 0.79 | 0.5 | 1 | 0.02 | -10.65 |
| ρ | 0.5 | 0.5 | 0 | 1 | 0.08 | -0.01 |
| Θ | 0.01 | 0.01 | -3.14 | 3.14 | 3.29 | -0.33 |
| **Gaussian** | P = 0.7914 | | | | | |
| SdM | 0.54 | 0.55 | 0 | 1 | 0.07 | -3.75 |
| MdL | 0.79 | 0.79 | 0.5 | 1 | 0.02 | -11.35 |
| ρ | 0.8 | 0.67 | 0 | 4.7 | 0.36 | 54.68 |
| Θ | 0.01 | 0.01 | -3.14 | 3.14 | 3.3 | -0.12 |
| **Ellipse 1:1** | P = 0.7213 | | | | | |
| SdM | 0.57 | 0.59 | 0 | 1 | 0.07 | -9.26 |
| MdL | 0.8 | 0.81 | 0.5 | 1 | 0.02 | -16.87 |
| ρ | 0.67 | 0.71 | 0 | 1 | 0.06 | -31.16 |
| Θ | 0 | 0.01 | -3.14 | 3.14 | 3.29 | -0.02 |
| **Ellipse 1:2** | P = 0.7929 | | | | | |
| SdM | 0.55 | 0.55 | 0 | 1 | 0.07 | -2.86 |
| MdL | 0.78 | 0.79 | 0.5 | 1 | 0.02 | -10.78 |
| ρ | 1.03 | 1 | 0 | 2 | 0.19 | 5 |
| Θ | 0 | -0.01 | -3.14 | 3.14 | 2.98 | 0.14 |

| Method \ Metrics | μ | Md | Min | Max | σ² | Skew |
|---|---|---|---|---|---|---|
| **Ellipse 1:3** | P = 0.863 | | | | | |
| SdM | 0.51 | 0.5 | 0 | 1 | 0.07 | 4.68 |
| MdL | 0.77 | 0.78 | 0.5 | 1 | 0.02 | -5.66 |
| ρ | 1.42 | 1.33 | 0.01 | 3 | 0.49 | 15.19 |
| Θ | 0 | -0.01 | -3.14 | 3.14 | 2.84 | 0.09 |
| **Rectangle 1:1** | P = 0.7257 | | | | | |
| SdM | 0.57 | 0.59 | 0 | 1 | 0.07 | -9.69 |
| MdL | 0.8 | 0.81 | 0.5 | 1 | 0.02 | -16.17 |
| ρ | 0.38 | 0.4 | 0 | 0.71 | 0.02 | -16.91 |
| Θ | 0 | -0.01 | -3.14 | 3.14 | 3.26 | 0.1 |
| **Rectangle 1:2** | P = 0.7981 | | | | | |
| SdM | 0.54 | 0.54 | 0 | 1 | 0.07 | -1.02 |
| MdL | 0.79 | 0.79 | 0.5 | 1 | 0.02 | -11.69 |
| ρ | 0.59 | 0.58 | 0 | 1.12 | 0.07 | -1.7 |
| Θ | 0.01 | 0.01 | -3.14 | 3.14 | 3.66 | -0.29 |
| **Rectangle 1:3** | P = 0.8672 | | | | | |
| SdM | 0.5 | 0.48 | 0 | 1 | 0.07 | 7.95 |
| MdL | 0.78 | 0.78 | 0.5 | 1 | 0.02 | -7.27 |
| ρ | 0.82 | 0.81 | 0.01 | 1.58 | 0.16 | 3.27 |
| Θ | -0.01 | 0 | -3.14 | 3.14 | 3.88 | 0.13 |
| **Fractal number** | P = 0.747 | | | | | |
| SdM | 0.58 | 0.6 | 0 | 1 | 0.06 | -10.94 |
| MdL | 0.78 | 0.79 | 0.5 | 1 | 0.02 | -10.1 |
| ρ | 25.56 | 24.19 | 0.04 | 93 | 171 | 30.35 |
| Θ | 0 | 0 | -3.14 | 3.14 | 3.29 | -0.06 |
| **Random quotient** | P = 0.7508 | | | | | |
| SdM | 0.58 | 0.59 | 0 | 1 | 0.06 | -9.38 |
| MdL | 0.79 | 0.79 | 0.5 | 1 | 0.02 | -11.2 |
| ρ | 0.9 | 0.82 | 0 | 3.3 | 0.24 | 68.95 |
| Θ | 0 | 0 | -3.14 | 3.14 | 3.29 | -0.13 |
| **L method** | P = 0.6401 | | | | | |
| SdM | 0.6 | 0.62 | 0 | 1 | 0.07 | -14.65 |
| MdL | 0.81 | 0.83 | 0.5 | 1 | 0.02 | -24.34 |
| ρ | 0.81 | 0.83 | 0.5 | 1 | 0.02 | -24.34 |
| Θ | 0.44 | 0.42 | 0 | 1.05 | 0.07 | 12.2 |
| **M method** | P = 0.8216 | | | | | |
| SdM | 0.68 | 0.73 | 0 | 1 | 0.05 | -34.69 |
| MdL | 0.72 | 0.7 | 0.5 | 1 | 0.02 | 19.53 |
| ρ | 0.86 | 0.87 | 0 | 1.73 | 0.13 | 0.54 |
| Θ | 0.62 | 0.62 | 0 | 1.56 | 0.1 | 4.58 |

The calculations in the Table are made for 175000 randomly created triangles.

This is the description of columns in the above Table:

1. *Method* - Randomization method, definition of sample space (in Bold font).
2. In the same Method column measured variables are shown (in Regular font) for different methods of randomization (sample spaces): ρ - radius and θ angle of a a polar coordinates of triangle vertices; SdM, and MdL are ratios of small to medium and medium to large side lengths of triangles respectively.
3. *Metrics*: µ - mean, Md - median, Min - minimum, Max - maximum, σ² - variance, Skew - skewness

These are different randomization methods:

1. *Generated* - triangles are generated by selecting big ∠A next value in [π/3, π] range by incrementing it by (2π/3)/<count>, where <count> is the number of random triangles. The remaining (π - ∠A) is randomly split into 2 more angle values with a restriction that neither of them are bigger than ∠A. Considering that the opposite to ∠A side a = 1 (or randomly selected), the lengths of other 2 sides are calculated. Coordinate for ∠A are randomly selected in a unit length circle. Considering the bisector of ∠A parallel to x-axis, the whole triangle is rotated clockwise at randomly selected angle in [0, 2π] range.
2. *Betrand #2* - polar coordinates: uniform distribution of the angle θ in [0, 2π] range and radius ρ in [0, 1] range
3. *Gaussian* - polar coordinates: uniform distribution of the angle θ in [0, 2π] range and normally distributed radius ρ in [0, 1] range for each vertex
4. *Ellipse* of a, b (1:1, 1:2, 1:3) radii - square roots of uniformly distributed lengths of minor **a** and major **b** axis, in case of **a** == **b** it translates to square root of uniform distribution of radius in [0, 1] range of circle and the angle θ in [0, 2π] range in polar coordinates (Betrand #3 point selection).
5. *Rectangle* a,b (1:1, 1:2, 1:3) - cartesian coordinates selected for X in range [-**b**/2; **b**/2] and for Y [-**a**/2; **a**/2] range (when **a**==**b** it is a square).
6. *Fractal* d (8) - the coordinate for each level from d to 1 is calculated as the sum of products **d** and a random real in the [-1,1] range. **d**=8 means that the plane is 8 by 8 embedded squares with side 2 (the coordinate is a random real in [-1,1] range)
7. *Quotient* n, x (0.09, 1.4) - calculates a radius **ρ** by dividing a random number in [0,1] range by a random number in [n,x] range, while angle θ is uniform in [0, 2π] range
8. L method - selects coordinates of large angle vertex per description in the L method section. Polar coordinates are randomly created in the x=[0.5, 1] and y=[0, sin(π/3)] rectangle area. If the radius ρ > 1 the selection is rejected and new random coordinates are generated to guarantee that the selected point is inside the ADE figure (see Diagram 2).
9. M method - selects coordinates of medium angle vertex per description in the M method section. Random angle θ in [0, π/3] range and ρ in [1, (1.0 + cos(2 * θ)) / cos(θ)) - points inside BDF figure on Diagram 3

*Acknowledgments*

The author would like to acknowledge and thank S. Minovitsky for illuminating and instructive discussions.

*References*

[Aer2014] Aerts, D., de Bianchi M.S., (2014), Solving the hard problem of Bertrand's paradox. arXiv:1403.4139v2 [physics.hist-ph] 27 Jun 2014
[Ban2009] N. Banerjee. (2009), Random Obtuse Triangles and Convex Quadrilaterals. Master of Science in Computation for Design and Optimization thesis. Massachusetts Institute of Technology 2009.
[Ber1889] Bertrand, J. (1889), Calcul des probabilités. Gauthier-Villars, Paris.
[Cha2006] G. Chaitin. (2006), Meta Math!, Vintage Books, NY, 2006.

[Che2023] R.A. Chechile. (2023), Bertrand's Paradox Resolution and Its Implications for the Bing–Fisher Problem. Mathematics 2023, 11, 3282

[DiP2010] Di Porto, P., Crosignani, B., Ciattoni, A., Liu, H. C. (2010), Bertrand's paradox: a physical solution. arXiv:1008.1878v1 [physics.data-an] 11 Aug 2010.

[Dod1893] C.L. Dodgson (Lewis Carroll), (1893) Curiosa Mathematica. Part II: Pillow-Problems, thought out during Sleepless Nights, Macmillan, L.

[Ede2015] A. Edelman, G. Strang, (2015) Random Triangle Theory with Geometry and Applications, arXiv:1501.03053v1 [math.HO] 9 Jan 2015.

[Gar2001] Gardner, M. (2001), The Colossal Book of Mathematics. W.W. Norton & Company, NY,.

[Ham1991] R.W. Hamming, (1991) The Art of Probability. For scientists and engineers. CRC Press 2018

[Hay2024] Hayrapetyan, A. (2024), A Note on Bertrand's paradox. academia.edu/118451522

[Jay1973] Jaynes, E. T. (1973), The Well-Posed Problem, Foundations of Physics, 3 (4): 477–493

[Lac1962] Lacey, O. L. (1962). The human organism as a random mechanism. Journal of General Psychology, 66, 321–325.

[Mar1994] Marinoff, L. (1994), A Resolution of Bertrand's Paradox. Source: Philosophy of Science, Vol. 61, No. 1, pp. 1-24. The University of Chicago Press on behalf of the Philosophy of Science Association

[Pet2019] Petroni, N.C., (2019), Thou shalt not say "at random" in vain: Bertrand's paradox exposed arXiv:1810.07805v3 [math.HO] 22 Aug 2019

[Por1994] S.Portnoy, (1994), A Lewis Carroll pillow problem: Probability of an obtuse triangle. Statistical Science, v.9, No. 2, 279-284, 1994.

[Wan2011] Wang, J., Jackson, R.. (2011) Resolving Bertrand's Probability Paradox. Int. J. Open Problems Compt. Math., Vol. 4, No. 3, September 2011.